\documentclass[a4paper]{amsart}

\RequirePackage{amsmath}
\RequirePackage{bm}
\RequirePackage{amssymb}
\RequirePackage{upref}
\RequirePackage{amsthm}
\RequirePackage{enumerate}
\RequirePackage{pb-diagram}
\RequirePackage{amsfonts}
\RequirePackage[mathscr]{eucal}
\RequirePackage{verbatim}
\RequirePackage{xr}
\RequirePackage{graphicx}
\usepackage{calc}
\usepackage{xspace}
\RequirePackage{color}
\RequirePackage{ifthen}
\RequirePackage{fancybox}

\newcommand{\cf}{cf.\@\xspace}

\newcommand{\al}{\alpha}

\newcommand{\de}{\delta }

\newcommand{\h}{\eta}

\newcommand{\ka}{\kappa}

\newcommand{\s}{\sigma}

\newcommand{\C}{\varGamma}

\newcommand{\F}{\varPhi}

\newcommand{\fv}[2]{#1\hspace{0pt}_{|_{#2}}}

\newcommand{\const}{\tup{const}}

\newcommand{\msp[1]}[1]{\mspace{#1mu}}

\newcommand{\R}[1][n+1]{{\protect\mathbb R}^{#1}}

\newcommand{\Hh}[1][n+1]{{\protect\mathbb H}^{#1}}
\newcommand{\Ss}[1][n+1]{{\protect\mathbb S}^{#1}}

\newcommand{\N}{{\protect\mathbb N}}

\newcommand{\eR}{\stackrel{\lower1ex \hbox{\rule{6.5pt}{0.5pt}}}{\msp[3]\R[]}}
\newcommand{\eN}{\stackrel{\lower1ex \hbox{\rule{6.5pt}{0.5pt}}}{\msp[1]\N}}
\newcommand{\eO}{\stackrel{\lower1ex \hbox{\rule{6pt}{0.5pt}}}{\msc O}}

\DeclareMathOperator{\graph}{graph}

\DeclareMathOperator{\osc}{osc}

\newcommand\ra{\rightarrow}

\newcommand\hra{\hookrightarrow}

\newcommand\pde[2]{\frac {\partial#1}{\partial#2}}
\newcommand\pd[3]{\frac {\partial#1}{\partial#2^#3}}   
 
\newcommand{\un}{\infty}
\newcommand{\A}{\forall}

\newcommand{\set}[2]{\{\,#1\colon #2\,\}}
\newcommand{\uu}{\cup}
\newcommand{\ii}{\cap}
\newcommand{\uuu}{\bigcup}

\newcommand{\uud}{ \stackrel{\lower 1ex \hbox {.}}{\uu}}
\newcommand{\uuud}[1]{ \stackrel{\lower 1ex \hbox {.}}{\uuu_{#1}}}
\newcommand\su{\subset}

\newcommand{\sminus}[1][28]{\raise 0.#1ex\hbox{$\scriptstyle\setminus$}}

\newcommand{\wed}{\wedge}

\newcommand{\abs}[1]{\lvert#1\rvert}

\newcommand{\spd}[2]{\protect\langle #1,#2\protect\rangle}

\newcommand{\tit}{\textit}

\newcommand{\tup}{\textup}

\newcommand{\mc}{\protect\mathcal}
\newcommand{\msc}{\protect\mathscr}

\providecommand{\bysame}{\makebox[3em]{\hrulefill}\thinspace}

\newcommand{\cq}[1]{\glqq{#1}\grqq\,}

\newcommand{\bt}{\begin{thm}}
\newcommand{\bl}{\begin{lem}}
\newcommand{\bc}{\begin{cor}}
\newcommand{\bd}{\begin{definition}}
\newcommand{\bpp}{\begin{prop}}
\newcommand{\br}{\begin{rem}}
\newcommand{\bn}{\begin{note}}
\newcommand{\be}{\begin{ex}}
\newcommand{\bes}{\begin{exs}}
\newcommand{\bb}{\begin{example}}
\newcommand{\bbs}{\begin{examples}}
\newcommand{\ba}{\begin{axiom}}
\newcommand{\bas}{\begin{assumption}}

\newcommand{\et}{\end{thm}}
\newcommand{\el}{\end{lem}}
\newcommand{\ec}{\end{cor}}
\newcommand{\ed}{\end{definition}}
\newcommand{\epp}{\end{prop}}
\newcommand{\er}{\end{rem}}
\newcommand{\en}{\end{note}}
\newcommand{\ee}{\end{ex}}
\newcommand{\ees}{\end{exs}}
\newcommand{\eb}{\end{example}}
\newcommand{\ebs}{\end{examples}}
\newcommand{\ea}{\end{axiom}}
\newcommand{\eas}{\end{assumption}}

\newcommand{\bp}{\begin{proof}}
\newcommand{\ep}{\end{proof}}
\newcommand{\eps}{\renewcommand{\qed}{}\end{proof}}

\newcommand{\bal}{\begin{align}}

\newcommand{\bi}[1][1.]{\begin{enumerate}[\upshape #1]}
\newcommand{\bia}[1][(1)]{\begin{enumerate}[\upshape #1]}
\newcommand{\bin}[1][1]{\begin{enumerate}[\upshape\bfseries #1]}
\newcommand{\bir}[1][(i)]{\begin{enumerate}[\upshape #1]}
\newcommand{\bic}[1][(i)]{\begin{enumerate}[\upshape\hspace{2\cma}#1]}
\newcommand{\bis}[2][1.]{\begin{enumerate}[\upshape\hspace{#2\parindent}#1]}
\newcommand{\ei}{\end{enumerate}}

\newcommand\ndots{\raise 0.47ex \hbox {,}\hskip0.06em\cdots %
     \raise 0.47ex \hbox {,}\hskip0.06em} 

\newcommand{\q}{\quad}
\newcommand{\qq}{\qquad}

\newcommand{\hp}{\hphantom}

\newcommand\nd{\noindent}

\newskip\Csmallskipamount                                                
\Csmallskipamount=\smallskipamount
\newskip\Cmedskipamount
\Cmedskipamount=\medskipamount
\newskip\Cbigskipamount
\Cbigskipamount=\bigskipamount

\newcommand\cvs{\vspace\Csmallskipamount}   
\newcommand\cvm{\vspace\Cmedskipamount}

\newskip\csa
\csa=\smallskipamount

\newskip\cma
\cma=\medskipamount

\newskip\cba
\cba=\bigskipamount

\newdimen\spt
\newcommand\citem{\cvs\advance\itemno by
1{(\romannumeral\the\itemno})\hskip3pt}
\newcommand{\bitem}{\cvm\nd\advance\itemno by
1{\bf\the\itemno}\hspace{\cma}}

\newcount\itemno
\newcommand{\las}[1]{\label{S:#1}}

\newcommand{\lae}[1]{\label{E:#1}}
\newcommand{\lat}[1]{\label{T:#1}}

\newcommand{\rs}[1]{Section~\ref{S:#1}}

\newcommand{\re}[1]{\eqref{E:#1}}

\newcommand{\fre}[1]{\eqref{E:#1} on page~\tup{\pageref{E:#1}}}

\newskip\thmskip
\thmskip=\parindent

\newskip\hsk
\newenvironment{hinw}{\labelsep=0pt\begin{list}{}{\labelsep=0pt\itemindent=0pt\labelwidth=0pt\leftmargin=\parindent\rightmargin=0pt\partopsep=\cba}%
\item\it\nopagebreak\nopagebreak}%
{\end{list}}

\newcommand\bh{\begin{hinw}}
\newcommand{\eh}{\end{hinw}}

\newtheoremstyle{normal}
  {\cba}
  {\cba}
  {}
  {\thmskip}
  {\bfseries}
  {.}
  {\hsk}
  {}

\newtheoremstyle{abschnitt}
  {\cba}
  {\cba}
  {}
  {\thmskip}
  {\bfseries}
  {.}
  {\hsk}
  {}

\newtheoremstyle{italic}
  {\cba}
  {\cba}
  {\itshape}
  {\thmskip}
  {\bfseries}
  {.}
  {\hsk}
  {}

\newtheoremstyle{aufgaben}
  {\cba}
  {\cba}
  {}
  {}
  {\normalsize\bfseries}
  {.}
  {\hsk}
  {}

\newtheoremstyle{break}
  {\cba}
  {\cba}
  {\itshape}
  {}
  {\bfseries}
  {.}
  {\newline}
  {}

\swapnumbers
\theoremstyle{italic}
\newtheorem{thm}[subsection]{Theorem}
\newtheorem{lem}[subsection]{Lemma}
\newtheorem{prop}[subsection]{Proposition}
\newtheorem{cor}[subsection]{Corollary}

\theoremstyle{normal}
\newtheorem{rem}[subsection]{Remark}
\newtheorem{definition}[subsection]{Definition}
\newtheorem{example}[subsection]{Example}
\newtheorem{examples}[subsection]{Examples}
\newtheorem{ex}[subsection]{Exercise}
\newtheorem{note}[subsection]{}
\newtheorem{axiom}[subsection]{Axiom}
\newtheorem{assumption}[subsection]{Assumption}

\theoremstyle{aufgaben}
\newtheorem{exs}[subsection]{Exercises}

\swapnumbers

\numberwithin{equation}{section}
\numberwithin{figure}{section}

\newenvironment{textequation}[1][0.8]
{\begin{equation}
\begin{aligned}
\begin{minipage}{#1\linewidth}}
{\end{minipage}
\end{aligned}
\end{equation}
\ignorespacesafterend}

\newcommand{\btext}{\begin{textequation}}
\newcommand{\etext}{\end{textequation}}

\def\hinweis{\@startsection{subsection}{2}%
 \z@{0.7\linespacing\@plus 0.5\linespacing}{0.7\linespacing}%
{\normalfont\itshape\indent}}

\newcounter{hours}\newcounter{minutes}
\newcommand{\printtime}{%
\setcounter{hours}{\time/60}%
\setcounter{minutes}{\time-\value{hours}*60}%
\ifthenelse{\value{minutes}<10}{\thehours :0\theminutes}{\thehours:\theminutes}}

\usepackage[german,english]{babel}
\usepackage{graphicx}
\RequirePackage{amsmath}
\RequirePackage{bm}
\RequirePackage{amssymb}
\RequirePackage{upref}
\RequirePackage{amsthm}
\RequirePackage{enumerate}
\RequirePackage{pb-diagram}
\RequirePackage{amsfonts}
\RequirePackage[mathscr]{eucal}
\RequirePackage{verbatim}
\RequirePackage{xr}
\RequirePackage{graphicx}
\usepackage{calc}
\usepackage{xspace}
\usepackage{dsfont}

\makeatletter
\RequirePackage{color}
\newcommand{\ann}[1]{\renewcommand{\@makefnmark}{\mbox{$^{\color{red}{\@thefnmark}}$}}%
\footnote {#1}}
\makeatother

\RequirePackage{upref}
\RequirePackage{amsthm}
\RequirePackage{enumerate}
\usepackage[mathscr]{eucal}

\usepackage{xr-hyper}

\usepackage{calc}

\newlength{\oddsidemarginlength}
\newlength{\topmarginlength}

\newcounter{numberoflines}
\newcounter{tempcc}
\oddsidemargin=\oddsidemarginlength
\evensidemargin=\oddsidemargin
\topmargin=\topmarginlength
\usepackage[colorlinks=true,linkcolor=blue,citecolor=blue,urlcolor=blue]{hyperref}  

\begin{document}

\flushbottom


\title[Pinching estimates]{Pinching estimates for dual flows in hyperbolic and de Sitter space}

\author{Claus Gerhardt}
\address{Ruprecht-Karls-Universit\"at, Institut f\"ur Angewandte Mathematik,
Im Neuenheimer Feld 294, 69120 Heidelberg, Germany}
\email{\href{mailto:gerhardt@math.uni-heidelberg.de}{gerhardt@math.uni-heidelberg.de}}
\urladdr{\href{http://www.math.uni-heidelberg.de/studinfo/gerhardt/}{http://www.math.uni-heidelberg.de/studinfo/gerhardt/}}
\thanks{This work was supported by the DFG}

%
\subjclass[2000]{35J60, 53C21, 53C44, 53C50, 58J05}
\keywords{dual flows, inverse  curvature flow, pinching estimates}
\date{\today}
%


\begin{abstract} 
We prove pinching estimates for dual flows provided the curvature function used in the inverse flow in de Sitter space is convex.
\end{abstract}

\maketitle

\tableofcontents

\setcounter{section}{0}
\section{Introduction}
In a recent paper \cite{cg:cfs} we considered dual flows in $\Ss$, where one flow is a direct flow
\begin{equation}
\dot x- -F\nu
\end{equation}
while the other is an inverse flow
\begin{equation}
\dot x=\tilde F^{-1}\nu,
\end{equation}
where $F$ is a smooth curvature function defined in $\C_+$ and $\tilde F$ its inverse
\begin{equation}
\tilde F(\ka_i)=\frac1{F(\ka_i^{-1})}.
\end{equation}
We assumed that $F$ is symmetric, monotone, positive, homogeneous of degree $1$, and that both $F$ and $\tilde F$ are concave; $F$ was also supposed to be strictly concave, \cf \cite[Definition 3.1]{cg:cfs} for a precise definition. Using the Gau{\ss} map in $\Ss$ which maps closed, strictly convex hypersurfaces $M$ to their polar sets $\tilde M$, we proved that the flow hypersurfaces of the above flows are polar sets of each other, if the initial hypersurfaces have this property, hence the name dual flows. Thus, if one problem is solvable for arbitrarily closed, strictly convex initial hypersurfaces the dual problem is also solvable. We were therefore able to prove that the direct flows contract to a round point and the inverse flows to an equator such that, after an appropriate rescaling, both flows converge to a geodesic sphere exponentially fast.

Let us emphasize that the result for the inverse flow could only be achieved because of the duality. We could tackle the inverse problem directly. However, Makowski and Scheuer \cite{scheuer:icfs} independently proved that the inverse flow converges to an equator in $C^{1,\al}$.

There also exist Gau{\ss} maps between hyperbolic space $\Hh$ and de Sitter space $N$ and vice versa mapping closed, strictly convex hypersurfaces of one space to closed, strictly convex hypersurfaces of the other space, \cf \rs{2} or \cite[Chapter 10.4]{cg:cp}. Moreover, the duality between direct flows in one space  and inverse flows in the other is also valid, \cf \cite{hao:dualflows}. When a direct flow in $\Hh$ with closed, strictly convex initial hypersurface $M$ contracts to a point $x_0$ the dual flow in $N$ will expand to a totally geodesic hypersurface which is isometric to $\Ss[n]$. After applying an isometry such that the point $x_0$ is the Beltrami point the hypersurfaces of the dual flow are all contained in $N_-$ and their limit hypersurface will be the slice
\begin{equation}
\{\tau=0\},
\end{equation}
compare \rs{2}.

As in the case of the dual flows in $\Ss$, the direct flow is easier to solve. A major ingredient in the proof is a pinching estimate for the principle curvatures. In \cite{hao:dualflows} pinching estimates could be derived provided the initial hypersurface in $\Hh$ is horoconvex, i.e., the principal curvatures satisfy the estimate
\begin{equation}\lae{1.5}
\ka_i\ge 1.
\end{equation}
We shall prove in this paper that pinching estimates can be derived provided the curvature function used in the inverse flow in $N$ is convex, monotone and positive without requiring that the principal curvatures of the initial hypersurface in $N$ satisfies the estimate
\begin{equation}
\ka_i\le 1
\end{equation}
which would be equivalent  to the relation \re{1.5} in $\Hh$. Hence, pinching estimates can be proved for merely strictly convex hypersurfaces.
\bt
Let $F\in C^\un(\C_+)\ii C^0(\bar\C_+)$ be monotone, positive, convex, symmetric and homogeneous of degree $1$, then the solution hypersurfaces of the flow
\begin{equation}
\dot x=-F^{-1}\nu
\end{equation}
in $N$ with strictly convex, spacelike, closed initial hypersurface $M_0$ satisfy the pinching estimate 
\begin{equation}
\frac{\ka_n}{\ka_1}\le c_0,
\end{equation}
where the $\ka_i$ are labelled such that
\begin{equation}
\ka_1\le\cdots\le\ka_n
\end{equation}
and where $c_0$ only depends on $M_0$ and $F(1,\ldots,1)$.
\et
In \cite{hao:dualflows} Yu has shown that the converge results described above are also valid for this configuration.

\section{Isometries in hyperbolic and de Sitter space and the Gau{\ss} maps}\las{2}
Hyperbolic and de Sitter space are \cq{spheres} in the Minkowski space $\R[n+1,1]$ defined by
\begin{equation}
\Hh=\set{x\in \R[n+1,1]}{\spd xx=-1\q\wed\q x^0>0}
\end{equation}
and
\begin{equation}
N=\set{x\in \R[n+1,1]}{\spd xx=1}.
\end{equation}
A point $x\in\R[n+1,1]$ has coordinates $x=(x^a)$, $0\le a\le n+1$, where $x^0$ is the time coordinate.

The isometries are the elements of the linear Lorentz group $\mc O(n+1,1)$ in $\R[n+1,1]$. For completeness we shall repeat the simple  proof that $\mc O(n+1,1)$ acts transitively in these spaces.
\bl
$\mc O(n+1,1)$ acts transitively in $\Hh$ and $N$.
\el
\bp
Let $x=(x^a)$ be an arbitrary point in $\Hh$. Using a rotation around the $x^0$-axis it can be mapped to the point $x_0=(x^0,x^1,0,\ldots,0)$. The map acting like
\begin{equation}
\begin{pmatrix}
\cosh\theta&\sinh\theta\\[\cma]
\sinh\theta&\cosh\theta
\end{pmatrix},\qq\theta\in\R[],
\end{equation}
in the coordinates $(x^0,x^1)$ and like the identity for the rest, is an isometry. Since 
\begin{equation}
x^0=\sqrt{1+\abs{x^1}^2}
\end{equation}
we can choose $\theta$ such that
\begin{equation}
(x^0,x^1,0,\dots,0)\q\ra \q p_0=(1,0,\ldots,0)
\end{equation}
which is known as the \tit{Beltrami point}.

By a similar proof, any point $x\in N$ can be mapped to $(0,1,0,\dots,0)$ by an isometry respecting the time orientation.
\ep
\bh
The Gau{\ss} maps
\eh

Let $M\su\Hh$ be a closed, strictly convex hypersurface, then $M$ can be considered to be a submanifold in $\R[n+1,1]$ of codimension $2$, i.e., 
\begin{equation}
x:M_0\ra \R[n+1,1].
\end{equation}

The Gaussian formula for $M$ then looks like
\begin{equation}\lae{9.3.3}
x_{ij}=g_{ij}x-h_{ij}\tilde x,
\end{equation}
where $g_{ij}$ is the induced metric, $h_{ij}$ the second fundamental form of $M$ considered as a hypersurface in $\Hh$, and $\tilde x$ is the representation of the (exterior\footnote{Notice that $M$ is orientable, since it is strictly convex, and  that for any closed, connected, orientable,  immersed hypersurface in $\Hh$ an exterior normal vector can be unambiguously defined.}) normal vector $\nu=(\nu^\al)$ of $M$ in $T(\Hh)$ as a vector in $T(\R[n+1,1])$.

The mapping
\begin{equation}\lae{2.8}
\tilde x:M_0\ra N
\end{equation}
is then the embedding of a strictly convex, closed, spacelike hypersurface $\tilde M$. We call this mapping the \tit{Gau{\ss} map}. In \cite[Theorem 10.4.4]{cg:cp} we have proved:
\bt\lat{9.4.4}
Let $x:M_0\ra M\su \Hh$ be a closed, connected, strictly convex hypersurface of class $C^m$, $m\ge 3$, then the Gau{\ss} map $\tilde x$ in \re{2.8} is the embedding of a closed, spacelike, achronal, strictly convex hypersurface $\tilde M\su N$ of class $C^{m-1}$.  

Viewing $\tilde M$ as a codimension $2$ submanifold in $\R[n+1,1]$, its Gaussian formula is
\begin{equation}\lae{9.4.25}
\tilde x_{ij}=-\tilde g_{ij}\tilde x+\tilde h_{ij} x,
\end{equation}
where $\tilde g_{ij}$, $\tilde h_{ij}$ are the metric and second fundamental form of the hypersurface $\tilde M\su N$, and $x=x(\xi)$ is the embedding of $M$ which also represents the future directed normal vector of $\tilde M$. The second fundamental form $\tilde h_{ij}$ is defined with respect to the future directed normal vector, where the time orientation of $N$ is inherited from $\R[n+1,1]$.

The second fundamental forms of $M$, $\tilde M$ and the corresponding principal curvatures $\ka_i$, $\tilde \ka_i$ satisfy
\begin{equation}\lae{9.4.26}
h_{ij}=\tilde h_{ij}=\spd{\tilde x_i}{x_j}
\end{equation}
and
\begin{equation}\lae{9.4.27}  
\tilde \ka_i=\ka_i^{-1}.
\end{equation}
\et
There also exists an inverse Gau{\ss} map mapping strictly convex hypersurfaces in $N$ to strictly convex hypersurfaces in $\Hh$:
\bt\lat{9.4.5}
Let $\tilde M\su N$ be a closed, connected, spacelike, strictly convex, embedded hypersurface of class $C^m$, $m\ge 3$, such that, when viewed as a codimension $2$ submanifold in $\R[n+1,1]$, its Gaussian formula is
\begin{equation}
\tilde x_{ij}=-\tilde g_{ij}\tilde x+\tilde h_{ij}x,
\end{equation}
where $\tilde x=\tilde x(\xi)$ is the embedding, $x$ the future directed\footnote{The time orientation is inherited from the ambient space $\R[n+1,1]$.} normal vector, and $\tilde g_{ij}$, $\tilde h_{ij}$ the induced metric and the second fundamental form of the hypersurface in $N$. Then we define the Gau{\ss} map as $x=x(\xi)$
\begin{equation}
x:\tilde M\ra \Hh\su \R[n+1,1].
\end{equation}
The Gau{\ss} map is the embedding of a closed, connected, strictly convex hypersurface $M$ in $\Hh$.

Let $g_{ij}$, $h_{ij}$ be the induced metric and second fundamental form of $M$, then, when viewed as a codimension $2$ submanifold, $M$ satisfies the relations
\begin{equation}\lae{9.4.37}
x_{ij}=g_{ij}x-h_{ij}\tilde x,
\end{equation}
\begin{equation}\lae{9.4.38}
h_{ij}=\tilde h_{ij}=\spd{x_i}{\tilde x_j},
\end{equation}
and
\begin{equation}\lae{9.4.39}
\ka _i=\tilde \ka_i^{-1},
\end{equation}
where $\ka_i$, $\tilde \ka_i$ are the corresponding principal curvatures.
\et
For a proof see \cite[Theorem 10.4.5]{cg:cp}.

When the Beltrami point is an interior point of the convex body defined by $M\su \Hh$, then the image if the Gau{\ss} map $\tilde x(M)$ is contained in
\begin{equation}
N_+=\set{x\in N}{\tau>0},
\end{equation}
where $\tau$ is the time coordinate in $N$ which is naturally defined by the embedding
\begin{equation}
N\hra \R[n+1,1].
\end{equation}
More precisely, we have proved in \cite[Theorem 10.4.9]{cg:cp}
\bt\lat{9.4.9}
The Gau{\ss} maps provide a bijective relation between the connected, closed, strictly convex hypersurfaces $M\su \Hh$ having the Beltrami point in the interior of their convex bodies and the spacelike, closed,  connected, strictly convex hypersurfaces $\tilde M\su N_+$. 

The geodesic spheres with center in the Beltrami point are mapped onto the coordinate slices $\{\tau=\const\}$. 
\et
\br
Our default definition for the second fundamental of a spacelike hypersurface in Lorentzian manifolds relies on the past directed normal vector and not on the future directed. Hence, we switch the light cone in $N$ such that the images of the Gau{\ss} map are strictly convex with respect to the past directed normal vector. But then, the coordinate used to prove the previous theorem is no longer future directed. In order to get a future directed coordinate system we have to replace the time function $\tau$ by $-\tau$. Then the image $\tilde x(M)$ is contained in $N_-$ instead of $N_+$; for details see the remarks at the end of \cite[page 307]{cg:cp}.
\er

\section{Convex curvature functions}
Let $F\in C^\un(\C_+)$ be monotone, convex, homogeneous of degree $1$ and normalized such that
\begin{equation}
F(1,\ldots,1)=1.
\end{equation}
Then the following estimates are valid
\begin{equation}
\frac1n H\le F\le \ka_n
\end{equation}
and
\begin{equation}
\sum_iF_i\le 1,
\end{equation}
where
\begin{equation}
F_i=\pd F{\ka}i
\end{equation}
and where the $\ka_i$ are labelled such that
\begin{equation}
\ka_1\le \cdots\le \ka_n.
\end{equation}
The $F_i$ satisfy the same ordering
\begin{equation}\lae{3.6} 
F_1\le\cdots\le F_n.
\end{equation}
Let $\tilde F$ be the inverse of $F$
\begin{equation}
\tilde F(\ka_i)=F^{-1}(\ka_i^{-1}),
\end{equation}
the $\tilde F$ is of class $(K)$ and hence strictly concave, \cf \cite[Lemma 2.2.12]{cg:cp} and \cite[Lemma 3.6]{cg:cfs}. 

From \re{3.6} we infer
\begin{equation}
F_i\ka_i\le F_j\ka_j,\q\text{if}\q \ka_i\le\ka_j,
\end{equation}
deducing further
\begin{equation}\lae{3.9}
\frac1n \ka_n\le F=\sum_iF_i\ka_i\le nF_n\ka_n,
\end{equation}
and therefore,
\begin{equation}\lae{3.10} 
\frac1{n^2}\le F_n\le 1.
\end{equation}
Define the relation
\begin{equation}
F\approx \ka_i
\end{equation}
if there are uniform constants $c_1,c_2$ such that
\begin{equation}
0<c_1\le \frac F{\ka_i}\le c_2,
\end{equation}
then  
\begin{equation}\lae{3.14}
F\approx \ka_n\q\wed\q \tilde F\approx \ka_1.
\end{equation}
Finally, we deduce form \re{3.9} and the homogeneity of $F$ that
\begin{equation}\lae{3.14.1}
1\le \sum_i\tilde F_i=F^{-2}\sum_iF_i\ka_i^2\le F^{-1}\ka_n\le n
\end{equation}
where the arguments of $\tilde F$ are $\ka_i^{-1}$.

Similarly we have 
\begin{equation}
1\ge \sum_i F_i=\tilde F^{-2}\tilde F_i\ka_i^2\ge \tilde F^{-1}\ka_1=F\ka_n^{-1}\ge \frac1n
\end{equation}

\section{Pinching estimates}
Let $N=N^{n+1}$ be the de Sitter space with metric given by
\begin{equation}
d\bar s^2=-d\tau^2+\cosh^2\tau\s_{ij}dx^idx^j,
\end{equation}
where $\tau\in\R[]$ and $\s_{ij}$ is the standard metric of $\Ss[n]$. We also assume that $\tau$ is a future directed time function; occasionally we also write $x^0$ instead of $\tau$. $N$ is space of constant curvature $K_N=1$.

In the following we shall formulate equations by using the general symbol $K_N$ for the curvature because of  transparency having in mind that $K_N=1$.

We consider the inverse curvature flow
\begin{equation}
\dot x=-F^{-1}\nu\equiv \F\nu,
\end{equation}
where $\nu$ is the past directed normal and $F$ is a curvature function defined in $\C_+$. $\F=\F(r)$ is the real function
\begin{equation}
\F=-r^{-1},\qq r>0.
\end{equation}
We shall consider smooth curvature functions $F$ which are monotone, convex, and homogeneous of degree $1$. The initial hypersurface $M_0$ of the flow should be closed, spacelike and strictly convex. As we pointed out in \rs{2} we may assume, by applying an isometry, 
\begin{equation}
M_0\su\{\tau<0\}.
\end{equation}
Spacelike hypersurfaces can be written as graphs over $\Ss[n]$, hence the leaves $M(t)$ of the flow can be expressed as
\begin{equation}
M(t)=\graph \fv{u(t,\cdot)}{\Ss[n]}. 
\end{equation}
Let
\begin{equation}
v^2=1-\cosh^{-2}u\,\s^{ij}u_ju_i,
\end{equation}
the $u$ satisfies the scalar curvature equation
\begin{equation}\lae{4.7}
\dot u =\pde ut=-v\F>0,
\end{equation}
hence we conclude
\bl
During the evolution the leaves $M(t)$ stay in compact set, or more precisely,
\begin{equation}\lae{4.8}
\inf_{M_0}u\le u<\sinh\pi. 
\end{equation}
\el
\bp
We shall prove that
\begin{equation}\lae{4.9}
\inf_{M(t)}u<0\qq\A\, t\in [0,T^*),
\end{equation}
from which the result would follow, since
\begin{equation}
\abs{Du}^2=\s^{ij}u_iu_j<\cosh^2u
\end{equation}
and hence
\begin{equation}
\osc u(t,\cdot)<\sinh\pi.
\end{equation}
Indeed, assume
\begin{equation}
\inf_{M(t)}u\ge 0
\end{equation}
for some $t$, then there would exist $\xi\in M(t)$ such that the second fundamental form in that point would satisfy
\begin{equation}
0<h_{ij}\le \bar h_{ij}\le 0,
\end{equation}
where $\bar h_{ij}$ is the second form of the coordinate slice
\begin{equation}
\{\tau =\inf_{M(t)}u\},
\end{equation}
a contradiction.
\ep
\br
In \cite{hao:dualflows} it is proved that the flow will converge to an umbilical maximal hypersurface $M$. Since any closed umbilical maximal hypersurface in $N$ is isometric to
\begin{equation}
\{\tau=0\},
\end{equation}
an isometric image of the flow will converge to the maximal coordinate slice, and, in view of the scalar curvature equation \re{4.7}, we can then deduce that the whole transformed flow is contained in $N_-$.
\er
The following evolution equation is well known
 \begin{equation}\lae{4.16}
\begin{aligned}
{\F}^\prime-\dot\F F^{ij}\F_{ij}&=
-\dot \F
F^{ij}h_{ik}h_j^k \F
+ K_N\dot\F F^{ij}g_{ij}\F,
\end{aligned}
\end{equation}where
\begin{equation}
\F^{\prime}=\frac{d}{dt}\F
\end{equation}
and
\begin{equation}
\dot\F=\frac{d}{dr}\F(r),
\end{equation}
see e.g, \cite[Lemma 2.3.4]{cg:cp}.

From \re{4.16} we infer
\begin{equation}\lae{4.19}
\dot F-\dot\F F^{ij} F_{;ij}=\dot\F F^{ij}h_{ik}h^k_jF-2F^{-1}\dot\F F^{ij}F_iF_j-K_N\dot\F F^{ij}g_{ij}F.
\end{equation}

Let  $(\tilde h^{ij})=(h_{ij})^{-1}$, then the mixed tensor $\tilde h^i_j$ satisfies the evolution equation
\begin{equation}\lae{4.20}
\begin{aligned}
\dot{\tilde h}^i_j&-\dot\F F^{kl}\tilde h^i_{j;kl}=\dot\F F^{kl}h_{kr}h^r_l\tilde h^i_j-\{\dot\F F -\F\}\de^i_j\\
&\hp =\;-K_N\{\dot\F F+\F\}\tilde h_{kj}\tilde h^{ki} +K_N\dot\F F^{kl}g_{kl}\tilde h^i_j\\
&\hp{=}\;-\{\dot\F F^{pq,kl}h_{pq;r}h_{kl;s}+2\dot\F F^{kl}\tilde h^{pq}h_{pk;r}h_{ql;s}+\Ddot\F F_rF_s\}\tilde h^{is}\tilde h^r_j,
\end{aligned}
\end{equation} 
\cf \cite[Lemma 2.4.3]{cg:cp}.

Finally, let $\h=\h(\tau)$ be a solution of the ordinary differential equation
\begin{equation}
\dot\h=-\frac{\bar H}n\h,
\end{equation}
where $\bar H$ is the mean curvature of the slices
\begin{equation}
\{x^0=\tau\},
\end{equation}
then the function
\begin{equation}
\chi=\tilde v\h(u),
\end{equation}
where $\tilde v-v^{-1}$,  satisfies
\begin{equation}\lae{4.24}
\begin{aligned}
\dot\chi-\dot\F F^{ij}\chi_{ij}&=-\dot\F F^{ij}h_{ik}h^k_j\chi+\{\dot\F F+\F\}\frac{\bar H}n\chi,
\end{aligned}
\end{equation}
\cf \cite[Lemma 10.5.5]{cg:cp}. In case of the de Sitter spacetime we can choose
\begin{equation}
\h=\cosh\tau,
\end{equation}
then
\begin{equation}
\dot\h=-\big(-\frac{\sinh\tau}{\cosh\tau}\big)\h.
\end{equation}
Applying the maximum principle we deduce from \re{4.8} and \re{4.24}
\bl
The quantities $\chi$ and $\tilde v$ are uniformly bounded during the evolution.
\el
From the equations \re{4.19} and \re{4.24} we infer
\bl
$F$ is uniformly bounded during the evolution.
\el
\bp
Let $0<T<T^*$ be arbitrary and apply the maximum principle to the function
\begin{equation}
w=\log F+\log\chi
\end{equation}
in
\begin{equation}
Q_T=[0,T]\times \Ss[n]
\end{equation}
assuming that there exist $(t_0,\xi_0)$, $0<t_0\le T$ such that
\begin{equation}
w(t_0,\xi_0)=\sup_{Q_T}w,
\end{equation}
then we have in $(t_0,\xi_0)$
\begin{equation}
0\le -K_N\dot \F F^{ij}g_{ij}<0,
\end{equation}
a contradiction; hence
\begin{equation}
w\le \sup_{M_0}w\q\text{in}\q Q_T
\end{equation}
from which the result immediately follows, since $\chi$ is uniformly positive
\begin{equation}
0<c_1\le \chi\le c_2<\un.
\end{equation}
\ep

We are now ready to prove the pinching estimates. First, let us observe that it suffices to prove
\begin{equation}
\frac F{\ka_1}\le \const,
\end{equation}
where the principal curvatures are labelled
\begin{equation}
\ka_1\le\cdots\le \ka_n,
\end{equation}
 in view of \fre{3.14}.
 \bt
 Let $F$ be monotone, convex, and homogeneous of degree $1$, then
 \begin{equation}
\ka_1^{-1} F\le \const
\end{equation}
during the evolution.
 \et
\bp
We shall show that the function
\begin{equation}
w=\log F+\log \tilde H-2nt,
\end{equation}
where
\begin{equation}
\tilde H=\tilde h^{ij}g_{ij},
\end{equation}
is uniformly bounded from above during the evolution by applying the maximum principle and using the fact that $T^*<\un$; the latter can easily be checked by comparing with the evolution of a coordinate slice, in view of \re{4.9}.

Let $0<T<T^*$ be arbitrary and suppose there exists $(t_0,\xi_0)$, $0<t_0\le T$, such that
\begin{equation}
w(t_0,\xi_0)=\sup_{QT}w.
\end{equation}
Applying the maximum principle we obtain
\begin{equation}
\begin{aligned}
0&\le 2F^{-2}F^{ij}h_{ik}h^k_j-2n F^{-1} -2n\\
&<2n -2n=0,
\end{aligned}
\end{equation}
in view of \fre{3.14.1}, a contradiction, proving the theorem. Here, we also used that the last term on the right-hand side of equation \re{4.20} is non-positive, \cf \cite[equ. (2.2.43) and Lemma 2.2.14]{cg:cp}.
\ep

\bibliographystyle{hamsplain}
\providecommand{\bysame}{\leavevmode\hbox to3em{\hrulefill}\thinspace}
\providecommand{\href}[2]{#2}



\end{document}